\documentclass[11pt]{article}
\usepackage{graphicx,array,amssymb,psfrag,fullpage}
\usepackage[round]{natbib}

\usepackage{algorithmic,algorithm}

\newcommand{\BEAS}{\begin{eqnarray*}}
\newcommand{\EEAS}{\end{eqnarray*}}
\newcommand{\BEA}{\begin{eqnarray}}
\newcommand{\EEA}{\end{eqnarray}}
\newcommand{\BEQ}{\begin{equation}}
\newcommand{\EEQ}{\end{equation}}
\newcommand{\BIT}{\begin{itemize}}
\newcommand{\EIT}{\end{itemize}}
\newcommand{\BNUM}{\begin{enumerate}}
\newcommand{\ENUM}{\end{enumerate}}

\newcommand{\BA}{\begin{array}}
\newcommand{\EA}{\end{array}}
\newcommand{\BC}{\begin{center}}
\newcommand{\EC}{\end{center}}

\newcommand{\ones}{\mathbf 1}

\newcommand{\reals}{{\mbox{\bf R}}}

\newcommand{\symm}{{\mbox{\bf S}}}  

\newcommand{\Card}{\mathop{\bf Card}}
\newcommand{\Tr}{\mathop{\bf Tr}}

\newcommand{\QED}{~~\rule[-1pt]{6pt}{6pt}}

\newcommand{\vect}{\mathop{\bf vec}}

\newtheorem{remark}[theorem]{Remark}

\newcounter{exno}

{\begin{quote}}{\end{quote}}

\makeatletter
\long\def\@makecaption#1#2{
   \vskip 9pt 
   \begin{small}
   \setbox\@tempboxa\hbox{{\bf #1:} #2}
   \ifdim \wd\@tempboxa > 5.5in
        \begin{center}
        \begin{minipage}[t]{5.5in}
        \addtolength{\baselineskip}{-0.95pt}
        {\bf #1:} #2 \par
        \addtolength{\baselineskip}{0.95pt}
        \end{minipage}
        \end{center}
   \else 
    \hbox to\hsize{\hfil\box\@tempboxa\hfil}  
   \fi
   \end{small}\par
}
\makeatother

\newcounter{oursection}

\newcounter{lecture}

\begin{document}

\title{A Pathwise Algorithm for Covariance Selection}

\author{
Vijay Krishnamurthy\thanks{ORFE, Princeton University, Princeton, NJ 08544, \texttt{kvijay@princeton.edu}}, Selin Damla Ahipa{\c s}ao{\u g}lu\thanks{ORFE, Princeton University, Princeton, NJ 08544, \texttt{sahipasa@princeton.edu}}, 
Alexandre d'Aspremont\thanks{ORFE, Princeton University, Princeton, NJ 08544, \texttt{aspremon@princeton.edu}}}

\maketitle

\begin{abstract}
Covariance selection seeks to estimate a covariance matrix by maximum likelihood while restricting the number of nonzero inverse covariance matrix coefficients. A single penalty parameter usually controls the tradeoff between log likelihood and sparsity in the inverse matrix. We describe an efficient algorithm for computing a full regularization path of solutions to this problem.\end{abstract}

\section{Introduction}
We consider the problem of estimating a covariance matrix from sample multivariate data by maximizing its likelihood, while penalizing the inverse covariance so that its graph is {\em sparse}. This problem is known as covariance selection and can be traced back at least to \citet{Demp72}. The coefficients of the inverse covariance matrix define the representation of a particular Gaussian distribution as a member of the exponential family, hence sparse maximum likelihood estimates of the inverse covariance yield sparse representations of the model in this class. Furthermore, in a Gaussian model, zeros in the inverse covariance matrix correspond to {\em conditionally} independent variables, so this penalized maximum likelihood procedure simultaneously stabilizes estimation and isolates {\em structure} in the underlying graphical model (see \citet{Laur96}).

Given a sample covariance matrix $\Sigma \in \symm_n$, the covariance selection problem is written as follows
\[\BA{ll}
\mbox{maximize} &  \log \det X -  \Tr (\Sigma X) - \rho \Card(X)  
\EA\]
in the matrix variable $X\in\symm_n$, where $\rho>0$ is a penalty parameter controlling sparsity and $\Card(X)$ is the number of nonzero elements in $X$. This is a combinatorially hard (non-convex) problem and, as in \citet{Dahl08,Bane06,Dahl05a}, we form the following convex relaxation
\BEQ\label{eq:covsel-intro}
\BA{ll}
\mbox{maximize}  &  \log \det X -  \Tr (\Sigma X) - \rho \|X\|_1 
\EA\EEQ
which is a convex problem in the matrix variable $X\in\symm_n$, where $\|X\|_1$ is the sum of absolute values of the coefficients of $X$ here. After scaling, the $\|X\|_1$ penalty can be understood as a convex lower bound on $\Card(X)$. Another completely different approach derived in \citet{Mein06} reconciles the local dependence structure inferred from $n$ distinct $\ell_1$-penalized regressions of a single variable against all the others. Both this approach and the convex relaxation (\ref{eq:covsel-intro}) have been shown to be consistent in \citet{Mein06} and \citet{Bane08} respectively.

In practice however, both methods are computationally challenging when $n$ gets large. Various algorithms have been employed to solve (\ref{eq:covsel-intro}) with~\citet{Dahl05a} using a custom interior point method and~\citet{Bane08} using a block coordinate descent method where each iteration required solving a LASSO-like problem, among others. This last method is efficiently implemented in the GLASSO package by \citet{Frie07} using coordinate descent algorithms from \citet{Frie07b} to solve the inner regression problems.

One key issue in all these methods is that there is no a priori obvious choice for the penalty parameter. In practice, at least a partial regularization path of solutions has to be computed, and this procedure is then repeated many times to get confidence bounds on the graph structure by cross-validation. Pathwise LASSO algorithms such as LARS by \citep{Efro04} can be used to get a full regularization path of solution using the method in \citet{Mein06} but this still requires solving and reconciling $n$ regularization paths on regression problems of dimension~$n$. 

Our contribution here is to formulate a pathwise algorithm for solving problem (\ref{eq:covsel-intro}) using numerical continuation methods (see \cite{Bach05} for an application in kernel learning). Each iteration requires solving a large structured linear system (predictor step) then improving precision using a block coordinate descent method (corrector step). Overall, the cost of moving from one solution to problem (\ref{eq:covsel-intro}) to another is typically much lower than that of solving two separate instances of (\ref{eq:covsel-intro}). We also derive a coordinate descent algorithm for solving the corrector step, where each iteration is closed-form and requires only solving a cubic equation. We illustrate the performance of our methods on several artificial and realistic data sets.

The paper is organized as follows. Section \ref{s:covsel} reviews some basic convex optimization results on the covariance selection problem in (\ref{eq:covsel-intro}). Our main pathwise algorithm is described in Section \ref{s:algo}. Finally, we present some numerical results in Section \ref{s:numres}.

\paragraph{Notation.} In what follows, we write $\symm_n$ for the set of symmetric matrices of dimension $n$. For a matrix $X\in\reals^{m\times n}$, we write $\|X\|_F$ its Frobenius norm, $\|X\|_1=\sum_{ij}|X_{ij}|$ the $\ell_1$ norm of its vector of coefficients, and $\Card(X)$ the number of nonzero coefficients in $X$.


\section{Covariance Selection}
\label{s:covsel} Starting from the convex relaxation defined above
\BEQ\label{eq:covsel}
\BA{ll}
\mbox{maximize}  &  \log \det X -  \Tr (\Sigma X) - \rho \|X\|_1 
\EA\EEQ
 in the variable  $X\in{\symm_n}$, where $\|X\|_1$ can be understood as a convex lower bound on the $\Card(X)$ function whenever $|X_{ij}|\leq 1$ (we can always scale $\rho$ otherwise). Let us write $X^\ast(\rho)$ for the optimal solution of problem (\ref{eq:covsel}). In what follows, we will seek to compute (or approximate) the entire regularization path of solutions $X^\ast(\rho)$, for $\rho \in \reals_+$. To remove the nonsmooth penalty, we can set $X=L-M$ and rewrite the problem above as 
\BEQ\label{eq:covsel-primal}
\BA{ll}
\mbox{maximize}  & \log \det (L -M) -  \Tr (\Sigma(L-M)) - \rho \ones^T(L +M)\ones\\
\mbox{subject to} & L_{ij},M_{ij} \geq 0,\quad i,j=1,\ldots,n,
\EA\EEQ
in the matrix variables $L, M \in \symm_n$. We can form the following dual to problem (\ref{eq:covsel}) as
\BEQ\label{eq:covsel-dual}
\BA{ll}
\mbox{minimize} &  -\log \det (U) - n \\
\mbox{subject to} & U_{ij} \leq \Sigma_{ij} + \rho, \quad i,j=1,\ldots,n, \\
                  & U_{ij} \geq \Sigma_{ij} - \rho, ,\quad i,j=1,\ldots,n,
\EA\EEQ
in the variable $U\in\symm_n$. As in \citet{Bach05} for example, in the spirit of barrier methods for interior point algorithms, we then form the following (unconstrained) regularized problem  
\BEQ\label{eq:regularized}
\min_{U\in\symm_n} -\log \det (U) -  t\left(\sum_{i,j=1}^n \log(\rho +\Sigma_{ij}
-U_{ij})  + \sum_{i,j=1}^n \log (\rho -\Sigma_{ij}+U_{ij})\right)
\EEQ
in the variable $U \in \symm_n$ and $t > 0$ specifies a desired tradeoff level between centrality (smoothness) and optimality. From every solution $U^\ast(t)$ corresponding to each $t>0$, the barrier formulation also produces an explicit \emph{dual} solution $(L^\ast(t), M^\ast(t))$ to Problem (\ref{eq:covsel-dual}). Indeed we can define matrices  $L, M \in\symm_n$ as follows
\[
L_{ij}(U,\rho) = \frac{t}{\rho+\Sigma_{ij}-U_{ij}} \quad \mbox{and} \quad 
M_{ij}(U,\rho) = \frac{t}{\rho-\Sigma_{ij}+U_{ij}}
\]
First order optimality conditions for problem (\ref{eq:regularized}) then imply 
\[
(L - M) = U^{-1}.
\] 
As $t$ tends to 0, problem (\ref{eq:regularized}) traces a central path towards the optimal solution to problem (\ref{eq:covsel-dual}). If we write $f(U)$ for the objective function of problem (\ref{eq:covsel-dual}) and call $p^\ast$ its optimal value, we get as in \cite[\S11.2.2]{Boyd03}
\[
f(U^\ast(t)) - p^\ast \leq 2n^2t
\]
hence $t$ can be understood as a surrogate duality gap when solving the dual problem (\ref{eq:covsel-dual}).


\section{Algorithm}
\label{s:algo} In this section we derive a Predictor-Corrector algorithm to approximate the entire path of solutions $X^\ast(\rho)$ when $\rho$ varies between 0 and $\max_{i}\Sigma_{ii}$ (beyond which the solution matrix is diagonal). Defining 
\[
H(U,\rho) = L(U,\rho) - M(U,\rho) - U^{-1}
\]
we trace the curve $H(U,\rho)=0$, the first order optimality condition for problem (\ref{eq:regularized}). Our pathwise covariance selection algorithm is defined in Algorithm \ref{alg:covpath}.

\begin{algorithm} 
\caption{Pathwise Covariance Selection} 
\label{alg:covpath} 
\begin{algorithmic} [1]
\REQUIRE $\Sigma\in\symm_m$
\STATE Start with $(U_0, \rho_0)$ s.t  $H(U_0,\rho_0) =0$.
\FOR{$i=1$ to $k$} 
\STATE\emph{Predictor Step}. Let $\rho_{i+1} = \rho_i + h$. Compute a tangent direction by solving the linear system
\[
\frac{\partial H}{\partial \rho}(U_i,\rho_i) + J(U_i,\rho_i) \frac{\partial U}{\partial\rho} = 0
\]
in $\partial U/\partial \rho\in\symm_n$, where $J(U_i,\rho_i)=\partial H(U,\rho)/\partial U\in\symm_{n^2}$ is the Jacobian matrix of the function $H(U,\rho)$. 
\STATE Update $U_{i+1} = U_i + h \partial U/\partial \rho$.
\STATE\emph{Corrector Step}. Solve problem (\ref{eq:regularized}) starting at $U = U_{i+1}$.
\ENDFOR 
\ENSURE Sequence of matrices $U_i$, $i=1,\ldots,k$.
\end{algorithmic} 
\end{algorithm} 

Typically in Algorithm \ref{alg:covpath}, $h$ is a small constant, $\rho_0=\max_{i}\Sigma_{ii}$, and $U_0$ is computed by solving a single (very sparse) instance of problem (\ref{eq:regularized}) for example.

\subsection{Predictor: conjugate Gradient method} \label{ss:pred}
In Algorithm \ref{alg:covpath}, the tangent direction in the predictor step is computed by solving a linear system  $Ax=b$ where $ A = (U^{-1} \otimes U^{-1} + D)$ and $D$ is a diagonal matrix. This system of equations has dimension $n^2$ and we solve it using the conjugate gradient (CG) method. 

\paragraph{CG iterations.} The most expensive operation in the CG iterations is the computation of a matrix vector product $Ap_k$, with $p_k \in \reals^{n^2 }$. Here however, we can exploit problem structure to compute this step efficiently. Observe that $(U^{-1} \otimes U^{-1})p_k = \vect(U^{-1} P_kU^{-1})$ when $p_k =\vect(P_k)$, so the  computation of the matrix vector product $Ap_k$ needs only $O(n^3)$ flops instead of $O(n^4)$. The CG method then needs at most $O(n^2)$ iterations to converge, leading to a total complexity of $O(n^5)$ for the predictor step. In practice, we will observe that CG needs considerably fewer iterations.

\paragraph{Stopping criterion.}
To speed up the computation of the predictor step, we can stop the conjugate gradient solver when the norm of the residual falls below the numerical tolerance $t$. In our experiments here, we stopped the solver after the residual decreases by two order of magnitudes.

\paragraph{Scaling \& warm start.}
Another option, much simpler than the predictor step detailed above, is warm starting. This means simply scaling the current solution to make it feasible for the problem after $\rho$ is updated. In practice, this method turns out to be as efficient as the predictor step as it allows us to follow the path starting from the sparse end (where more interesting solutions are located). Here, we start the algorithm from the sparsest possible solution, a diagonal matrix $U$ such that
\[
U_{ii}=\Sigma_{ii} +(1-\epsilon)\rho_\mathrm{max} I, \quad i=1,\ldots,n,
\]
where $\rho_\mathrm{max}=\max_i \Sigma_{ii}$. Suppose now that iteration $k$ of the algorithm produced a matrix solution $U_k$ corresponding to a penalty $\rho_k$, the algorithm with (lower) penalty $\rho_{k+1}$ is started at the matrix
\[
U=(1-\rho_{k+1}/\rho_{k})\Sigma+(\rho_{k+1}/\rho_{k})U_k
\]
which is a feasible starting point for the corrector problem that follows. This is the method that was implemented in the final version of our code and that is used in the numerical experiments detailed in the numerical section.

\subsection{Corrector: block coordinate descent} \label{ss:corr}
For small size problems, we can use Newton's method to solve problem~(\ref{eq:regularized}). However  from a computational perspective, this approach is not practical for large values of $n$. We can simplify iterations using a block coordinate descent algorithm that updates one row/column of the matrix in each iteration (\citet{Bane08}). Let us partition the matrices $U$ and $\Sigma$ as  
\[
U = \left(\BA{cc}V&u\\ 
u^T&w\EA\right)
\quad \mbox{and} \quad
S = \left(\BA{cc}A&b\\ 
b^T&c\EA\right)
\]
We keep $V$ fixed in each iteration and solve for $u$ and $w$. Without loss of generality, we can always assume that we are updating the last row/column.  
\paragraph{Algorithm.} Problem (\ref{eq:regularized}) can be written in block format as: 
\BEQ\label{eq:block-regularized}
\BA{ll}
\mbox{minimize} & -\log(w -u^TV^{-1}u)  -t(\log(\rho+c-w) +\log (\rho-c+w))\\
\\
& -  2t\left(\sum_{i} \log(\rho+b_i -u_i) + \sum_{i} \log(\rho -b_{i}+u_{i})\right)  
\EA
\EEQ
in the variables $u \in \reals^{(n-1)}$ and $w \in \reals$. Here  $V \in\symm^{(n-1)}$ is kept fixed in each iteration. 
\begin{algorithm} 
\caption{Block coordinate descent corrector steps} 
\label{alg:block-corr} 
\begin{algorithmic} [1]
\REQUIRE $U_0,~\Sigma\in\symm_n$
\FOR{$i=1$ to $k$} 
\STATE Pick the row and column to update. 
\STATE Solve the inner problem (\ref{eq:block-regularized}) using coordinate descent (each coordinate descent step requires solving a cubic equation).
\STATE Update $U^{-1}$.
\ENDFOR 
\ENSURE A matrix $U_k$ solving (\ref{eq:regularized}).
\end{algorithmic} 
\end{algorithm} 
We use the Sherman-Woodbury-Morrison (SWM) formula (see \citet[\S C.4.3]{Boyd03}) to efficiently update $U^{-1}$ at each iteration, so it suffices to compute the full inverse only once at the beginning of the path. The choice and order of row/column updates significantly affects performance. Although predicting the effect of a whole~$i^{th}$ row/column update is numerically expensive, we use the fact that the impact of updating diagonal coefficients usually dominates all others and can be computed explicitly at a very low computational cost. It corresponds to the maximum improvement in the dual objective function that can be achieved by updating the current solution $U$ to $U + w e_ie_i^T $, where $e_i$ is the~$i^{th}$ unit vector. The objective function value is a decreasing function of $w$ and $w$ must be lower than $\rho+\Sigma_{ii}-U_{ii}$ to preserve dual feasibility, so updating the $i^{th}$ diagonal coefficient will decrease the objective by $\delta_i=(\rho + \Sigma_{ii}-U_{ii})U^{-1}_{ii}$ after minimizing over $w$. In practice, updating the top 10\% row/columns with largest~$\delta$ is often enough to reach our precision target and very significantly speeds-up computations. We also solve the inner problem (\ref{eq:block-regularized}) by a coordinate descent method (as in \citet{Frie07b}), taking advantage of the fact that a point minimizing (\ref{eq:block-regularized}) over a single coordinate can be computed in closed-form by solving a cubic equation. Suppose $(u,w)$ is the current point and that we wish to optimize coordinate $u_j$ of the vector $u$, we define
\BEQ\label{eq:coefs}
\BA{ll} 
\alpha = - V^{-1}_{jj} \\
\beta = -2u_j(\sum_{k\neq j} V^{-1}_{kj}u_k) \\
\gamma =  w - u^TV^{-1}u - \alpha u_j - \beta u_j^2
\EA\EEQ
The optimality conditions imply that the the optimal $u_j^\ast$ must satisfy the following cubic equation
\BEQ\label{eq:cubic}
\BA{ll} 
p_1x^3 + p_2x^2 +p_3x +  p_4= 0
\EA\EEQ
where
\[\BA{ll} 
p_1 = 2(1+2t)\alpha,~p_2 = (1+4t)\beta - 4(1+2t)\alpha b_j \\
p_3 = 4t\gamma  - 2(1+2t)\beta b_j    + 2\alpha(b^2_j -2\rho^2),~p_4 = \beta(b^2_j-\rho^2) -4t\gamma b_j.
\EA\]
Similarly the diagonal update  $w$ satisfies the following quadratic equation. 
\[\BA{ll} 
(1+2t)w^2 -2( t(u^TV^{-1}u) +c(1+t))w + c^2 -\rho^2 + 2tc(u^TV^{-1}u) = 0
\EA\]
Here too, the order in which we optimize the coordinates has a significant impact.

\paragraph{Dual block problem.} We can derive a dual to problem (\ref{eq:block-regularized}) by rewriting it as a constrained optimization problem to get 
\BEQ\label{eq:dummy-block-regularized}
\BA{ll}
\mbox{minimize} & -\log x_1 -t(\log x_2 +\log x_3) -  2t\left(\sum_{i} (\log y_i + \log z_i) \right)\\
\mbox{subject to}  & x_1 \leq  w - u^TV^{-1}u \\
                               & x_2 = \rho + c -w,~x_3 = \rho -c +w \\
                               & y_i = \rho +b_i - u_i,~z_i = \rho -b_i +u_i
\EA
\EEQ
in the variables $u \in \reals^{(n-1)}, w \in \reals, x \in \reals^3, y \in \reals ^{(n-1)}, z \in \reals^{(n-1)}$. The dual to problem (\ref{eq:dummy-block-regularized}) is written
\BEQ\label{eq:dual-block-regularized}
\BA{ll}
\mbox{maximize} & 1 + 2t(2n-1)+ \log \alpha_1 - \alpha_2(\rho+c) - \alpha_3(\rho-c)\\ 
& - \sum_{i} \left ( \beta_i(\rho+b_i) + \eta_i(\rho-b_i)\right)\\
                    & + t\log(\alpha_2/t)+t\log(\alpha_3/t) + 2t\left(\sum_i\left(\log(\beta_i/2t) +    \log(\eta_i/2t) \right)\right)\\
\mbox{subject to}  &\alpha_1 = \alpha_2 - \alpha_3 \\
                               & \alpha_1 \geq 0      
\EA
\EEQ
in the variables  $\alpha \in \reals^3, \beta \in \reals ^{(n-1)}$ and $\eta \in \reals^{(n-1)}$. Surrogate dual points then produce an explicit stopping criterion. 

\subsection{Complexity}
Solving for the predictor step using conjugate gradient as in $\S\ref{ss:pred}$ requires $O(n^2)$ matrix products (at a cost of $O(n^3)$ each) in the worst-case, but the number of iterations necessary to get a good estimate of the predictor is typically much lower (cf. experiments in the next section). Scaling and warm start on the other hand has complexity $O(n^2)$. The inner and outer loops of the corrector step are solved using coordinate descent, with each coordinate iteration requiring the (explicit) solution of a cubic equation. 

Results on the convergence of the coordinate descent in the smooth case can be traced back at least to \citep{Luo92} or \citep{Tsen01}, who focus on local linear convergence in the strictly convex case. More precise convergence bounds have been derived in \cite{Nest10} who shows linear convergence (i.e. with complexity growing as $\log(1/\epsilon)$) of a randomized variant of coordinate descent for strongly convex functions, and a complexity bound growing proportionally to $1/\epsilon$ when the gradient is Lipschitz continuous coordinatewise. Unfortunately, because it uses a randomized step selection strategy, the algorithm in its standard form is inefficient in our case here, as it requires too many SWM matrix updates to switch between columns. Optimizing the algorithm in \cite{Nest10} to adapt it to our problem (e.g. by adjusting the variable selection probabilities to account for the relative cost of switching columns) is a potentially promising research direction. 

The complexity of our algorithm can be summarized as follows.
\begin{itemize}
\item Because our main objective function is strictly convex, our algorithm converges locally linearly, but we have no explicit bound on the total number of iterations required.
\item Starting the algorithm requires forming the inverse matrix $V^{-1}$ at a cost of $O(n^3)$.
\item Each iteration requires solving a cubic equation for each coordinatewise minimization problem to form the coefficients in (\ref{eq:coefs}), at a cost of $O(n^2)$. Updating the problem to switch from one iteration to the next using SWM updates then costs $O(n^2)$. This means that scanning the full matrix with coordinate descent requires $O(n^4)$ flops.
\end{itemize}
While the lack of precise complexity bound is a clear shortcoming of our choice of algorithm for solving the corrector step, as discussed by \cite{Nest08}, algorithm choices are usually guided by the type of operations (projections, barrier computations, inner optimization problems) that can be solved very efficiently or in closed-form. In our case here, it turns out that coordinate descent iterations can be performed very fast, in closed-form (by solving cubic equations), which seems to provide a clear (empirical) complexity advantage to this technique.


\section{Numerical Results} \label{s:numres} 
We compare the numerical performance of several methods for computing a full regularization path of solutions to problem (\ref{eq:covsel}) on several realistic data sets: the senator votes covariance matrix from \citet{Bane06}, the {\it Science} topic model in \citet{Blei07} with 50 topics, the covariance matrix of 20 foreign exchange rates, the UCI SPECTF heart dataset (diagnosing of cardiac images), the UCI LIBRAS hand movement dataset and the UCI HillValley dataset. We compute a path of solutions using the methods detailed here (Covpath) and repeat this experiment using the Glasso path code \citet{Frie07} which restarts the covariance selection problem at $\rho + \epsilon$ at the current solution of (\ref{eq:covsel}) obtained at $\rho$. We also tested the smooth first order code with warm-start ASPG described in \citep{Lu10a} as well as the greedy algorithm SINCO by \citet{Sche09}. Note that the later only identifies good sparsity patterns but does not (directly) produce feasible solutions to problem (\ref{eq:covsel-dual}). Our prototype code here is written in MATLAB (except for a few steps in~C), ASPG and SINCO are also written in MATLAB, while Glasso is compiled from Fortran and interfaced with R. We use the scaling/warm-start approach detailed in \S\ref{s:algo} and scan the full set of variables at each iteration of the block-coordinate descent algorithm (optimizing over the 10\% most promising variables sometimes significantly speeds up computations but is more unstable), so the results reported here describe the behavior of the most robust implementation of our algorithm. We report CPU time (in seconds) versus problem dimension in Table \ref{tab:cpu}.  Unfortunately, Glasso does not use the duality gap as a stopping criterion but rather lack of progress (average absolute parameter change less than $10^{-4}$). Glasso fails to converge on the HillValley example.

\begin{table}[H]
\BC
\begin{tabular}{r|rrrrr}
Dataset & Dimension & Covpath & Glasso & ASPG & SINCO\\
\hline
Interest Rates & 20 & 0.036 & 0.200 & 0.30 & \bf 0.007\\
FXData & 20 & \bf 0.016 & 1.467 & 4.88 & 0.109 \\
Heart & 44 & \bf 0.244 & 2.400 & 11.25 & 5.895\\
ScienceTopics & 50 & \bf 0.026 & 2.626 & 11.58 & 5.233\\
Libras & 91 & \bf 0.060 & 3.329 &  35.80 & 40.690\\
HillValley & 100 & \bf 0.068 & - & 47.22 & 68.815\\
Senator & 102 & \bf 4.003 & 5.208 & 10.44 & 5.092\\
\end{tabular}
\caption{{\small CPU time (in seconds) versus problem type for computing a regularization path for 50 values of the penalty $\rho$, using the path following method detailed here (Covpath), the Glasso code with warm-start (Glasso), the pathwise code (ASPG) in \citep{Lu10a} and the SINCO greedy code by \citet{Sche09}.} \label{tab:cpu}}
\EC
\end{table}

As in \cite{Bane08}, to test the behavior of the algorithm on examples with known graphs, we also sample sparse random matrices with Gaussian coefficients, add multiples of the identity to make them positive semidefinite, then use the inverse matrix as our sample matrix~$\Sigma$. We use these examples to study the performance of the various algorithms listed above on increasingly large problems. Computing times are listed in Table~\ref{tab:cpu10}, for a path of length 10, and Table~\ref{tab:cpu50} for a path of length 50. The penalty coefficients $\rho$ are chosen to produce a target sparsity around 10\%.

\begin{table}[H]
\BC
\begin{tabular}{r|rrrr}
Dimension & Covpath & Glasso & ASPG & SINCO\\
\hline
20&\bf0.0042&2.32&0.53&0.22\\
50&\bf0.0037&0.59&4.11&3.80\\
100&\bf0.0154&1.11&13.36&13.58\\
200&\bf0.0882&4.73&73.24&61.02\\
300&\bf0.2035&13.52&271.05&133.99\\
\end{tabular}
\caption{{\small CPU time (in seconds) versus problem dimension for computing a regularization path for 10 values of the penalty $\rho$, using the path following method detailed here (Covpath), the Glasso code with warm-start (Glasso), the pathwise code (ASPG) in \citep{Lu10a} and the SINCO greedy code by \citet{Sche09} on randomly generated problems.} \label{tab:cpu10}}
\EC
\end{table}

\begin{table}[H]
\BC
\begin{tabular}{r|rrrr}
Dimension & Covpath & Glasso & ASPG & SINCO\\
\hline
20&\bf0.0101&0.64&2.66&1.1827\\
50&\bf0.0491&1.91&23.2&22.0436\\
100&\bf0.0888&10.60&140.75&122.4048\\
200&\bf0.3195&61.46&681.72&451.6725\\
300&\bf0.8322&519.05&5203.46&1121.0408\\
\end{tabular}
\caption{{\small CPU time (in seconds) versus problem dimension for computing a regularization path for 50 values of the penalty $\rho$, using the path following method detailed here (Covpath), the Glasso code with warm-start (Glasso), the pathwise code (ASPG) in \citep{Lu10a} and the SINCO greedy code by \citet{Sche09}  on randomly generated problems.} \label{tab:cpu50}}
\EC
\end{table}

\begin{figure}[h]
\begin{center}
\includegraphics[width=\textwidth]{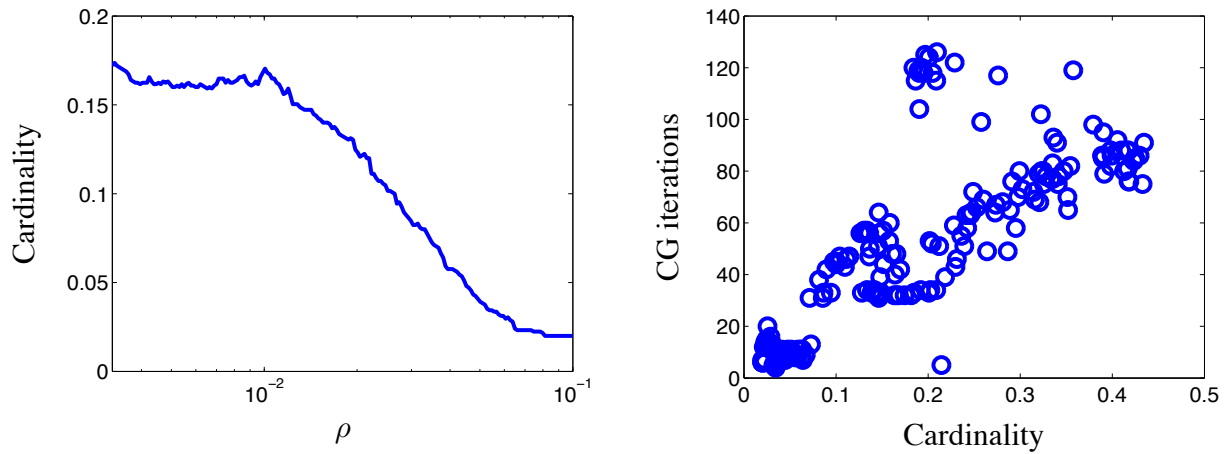} 
\caption{{\em Left:} We plot the fraction of nonzero coefficients in the inverse covariance versus penalty parameter $\rho$, along a path of solutions to problem (\ref{eq:covsel}). {\em Right:} Number of conjugate gradient iterations required to compute the predictor step versus number of nonzero coefficients in the inverse covariance matrix. \label{fig:nnz}}
\end{center}
\end{figure}

In Figure~\ref{fig:nnz}, we plot the number of nonzero coefficients (cardinality) in the inverse covariance versus penalty parameter~$\rho$, along a path of solutions to problem (\ref{eq:covsel}). We observe that the solution cardinality appears to be linear in the log of the regularization parameter. We then plot the number of conjugate gradient iterations required to compute the predictor in~\S\ref{ss:pred} versus number of nonzero coefficients in the inverse covariance matrix. We notice that the number of CG iterations decreases significantly for sparse matrices, which makes computing predictor directions faster at the sparse (i.e. interesting) end of the regularization path. Nevertheless, the complexity of corrector steps dominates the total complexity of the algorithm and there was little difference in computing time between using the scaling method detailed in \S\ref{s:algo} and using the predictor step, hence the final version of our code and the CPU time results listed here make use of scaling/warm-start exclusively, which is more robust.

\section{Online Covariance Selection}
\label{s:onl} 

In this section, we will briefly discuss the {\em online} version of the Covariance Selection problem. This version arises if we obtain a better estimate of the covariance matrix after the problem is already solved for a set of parameter values. We will assume that the new (positive definite) covariance matrix $\hat{\Sigma}$ is the sum of the old covariance matrix $\Sigma$ and an arbitrary symmetric matrix $C$. With such a change, the `new' dual problem can be written as 

\BEQ\label{eq:covsel-dual-onl}
\BA{ll}
\mbox{minimize} &  -\log \det (U) - n \\
\mbox{subject to} & U_{ij} \leq \rho +\Sigma_{ij} +\mu C_{ij}, \quad i,j=1,\ldots,n, \\
                  & U_{ij} \geq \Sigma_{ij} +\mu C_{ij} - \rho, ,\quad i,j=1,\ldots,n,
\EA\EEQ
in the variable $U\in\symm_n$, where $\rho$ is a parameter value for which the corresponding optimal solution is already calculated with the old covariance matrix $\Sigma$. The problem is parametrized with $\mu$, so that $\mu = 0$ gives the original problem whereas $\mu=1$ corresponds to the new problem. 

For many applications, one would expect $C$ to be small and the optimal solution $U^*$ of the original problem to be close to the optimal solution of the new problem, say $\hat{U}^*$. Hence, regardless of the algorithm, $U^*$ should be used as an initial solution instead of solving the problem from scratch.  

In the spirit of the barrier methods and the predictor-corrector method that we have devised in this chapter, we can develop a predictor-corrector algorithm to solve the online version of the problem fast as follows. We form a parametrized version of the regularized problem 

\BEQ\label{eq:regularized-ONL}
\BA{ll}
\displaystyle \min_{U\in\symm_n} & -\log \det (U) -  t\sum_{i,j=1}^n \log(\rho +\Sigma_{ij} +\mu C_{ij}
-U_{ij})\\
&    - t \sum_{i,j=1}^n \log (\rho -\Sigma_{ij} - \mu C_{ij} + U_{ij})
\EA\EEQ
in the variable $U \in \symm_n$ and $t > 0$ the tradeoff level as before. 
Let us define matrices  $\hat{L}, \hat{M} \in\symm_n$ as follows
\[
\hat{L}_{ij}(U,\mu) = \frac{t}{\rho+\Sigma_{ij} + \mu C_{ij}-U_{ij}} \quad \mbox{and} \quad 
\hat{M}_{ij}(U,\mu) = \frac{t}{\rho-\Sigma_{ij}-\mu C_{ij}+U_{ij}}
\]
As before, optimal $\hat{L}$ and $\hat{M}$ should satisfy $ (\hat{L} -\hat{ M}) = U^{-1}$, and problem (\ref{eq:regularized-ONL}) traces a central path towards the optimal solution to problem (\ref{eq:covsel-dual-onl}) as $t$ goes to 0.

Defining 
\[
\hat{H}(U,\mu) = \hat{L}(U,\mu) - \hat{M}(U,\mu) - U^{-1},
\]
we trace the curve $\hat{H}(U,\mu)=0$, the first order optimality condition for problem (\ref{eq:regularized-ONL}), from the solution for the original problem to one for the new problem as $\mu$ goes from 0 to 1. The resulting predictor-corrector algorithm is Algorithm \ref{alg:covpath-ONL}, which solves the online version efficiently.

\begin{algorithm}[ht]
\caption{Online Pathwise Covariance Selection} 
\label{alg:covpath-ONL} 
\begin{algorithmic} [1]
\REQUIRE $\Sigma, U^*\in\symm_m$, $\rho \in \reals$, and $c \in \reals^{n \times r}$. 
\STATE Start with $(U_0, \mu_0)$ s.t  $\hat{H}(U_0,\mu_0) =0$, specifically, set $\mu_0=0$ and $U_0 = U^*$.
\FOR{$i=1$ to $k$} 
\STATE\emph{Predictor Step}. Let $\mu_{i+1} = \mu_i + 1/k$. Compute a tangent direction by solving the linear system
\[
\frac{\partial \hat{H}}{\partial \mu}(U_i,\mu_i) + J(U_i,\mu_i) \frac{\partial U}{\partial\mu} = 0
\]
in $\partial U/\partial \mu\in\symm_n$, where $J(U_i,\mu_i)=\partial \hat{H}(U,\mu)/\partial U\in\symm_{n^2}$ is the Jacobian matrix of the function $\hat{H}(U,\mu)$. 
\STATE Update $U_{i+1} = U_i + (\partial U/\partial \mu)/k$.
\STATE\emph{Corrector Step}. Solve problem (\ref{eq:regularized-ONL}) for $\mu_{i+1}$ starting at $U = U_{i+1}$.
\ENDFOR 
\ENSURE Matrix $U_k$ that solves Problem (\ref{eq:covsel-dual-onl}).
\end{algorithmic} 
\end{algorithm} 

As for the offline version, the most demanding computation in this algorithm is the calculation of the tangent direction which can be carried out by the CG method discussed above. When carefully implemented and tuned, it produces a solution for the new problem very fast. Although one can try different values of $k$, setting $k=1$, and applying one step of the algorithm is usually enough in practice. This algorithm, and the online approach discussed in this section in general, would be especially useful and sometimes necessary for very large data sets as solving the problem from scratch is an expensive task for such problems and should be avoided whenever possible.

\section*{Acknowledgements}
The authors are grateful to two anonymous referees whose comments significantly improved the paper. The authors would also like to acknowledge support from NSF grants SES-0835550 (CDI), CMMI-0844795 (CAREER), CMMI-0968842, a Peek junior faculty fellowship, a Howard B. Wentz Jr. award and a gift from Google.

\bibliographystyle{plainnat}
\bibliography{MainPerso}

\end{document}